\documentclass[12pt]{amsart}  
\usepackage{graphicx}  
\usepackage{amsmath,amsthm}  
\usepackage{amssymb,amsmath} 
\usepackage{amscd}  
\usepackage{epsfig}  
\usepackage{graphicx}  
\usepackage{tabularx} 
\usepackage{color} 
\DeclareGraphicsExtensions{.jpg, .eps} 
\usepackage{textcomp} 
  
\newcommand{\C}{{\mathcal C}}  
  
\newcommand{\p}{{\mathcal P}}  
  
\newtheorem{theorem}{Theorem}  
\newtheorem{lemma}[theorem]{Lemma}  
\newtheorem{proposition}[theorem]{Proposition}  
\newtheorem{definition}{Definition}  
\newtheorem{corollary}[theorem]{Corollary}

\begin{document} 
 
\title{Rhombus filtrations and Rauzy algebras} 
 
 \author{Alex Clark}   
\address{University of North Texas,  
P.O.Box 311430, Denton, Texas 76203-1430, USA} 
 \email{alexc@unt.edu}

\author{Karin Erdmann}   
\address{Mathematical Institute, 24-29 St Giles',   
 Oxford OX1 3LB, UK}   
 \email{erdmann@maths.ox.ac.uk}

\author{Sibylle Schroll}   
\address{Mathematical Institute, 24-29 St. Giles',    
Oxford OX1 3LB, UK}   
\email{schroll@maths.ox.ac.uk}

\thanks{The third author acknowledges support through a Marie Curie   
Fellowship}   

\begin{abstract}   
Peach introduced rhombal algebras associated to quivers given by tilings  
of the plane by rhombi. We develop general techniques to analyze rhombal  
algebras, including a filtration by what we call rhombus modules.  
 
We introduce a way to relate the infinite-dimensional rhombal 
algebra corresponding to a complete tiling of the plane to finite-dimensional 
algebras corresponding to finite portions of the tiling. Throughout, we apply our general  
techniques to the special case of the Rauzy tiling,  which is built in stages reflecting 
an underlying self-similarity. Exploiting this self-similar structure allows us 
to uncover interesting features of the associated finite-dimensional algebras, including some of the  
tree classes in the stable Auslander-Reiten quiver. 
\end{abstract}

\maketitle

\parindent0pt  
\baselineskip=18pt  
\section{Introduction}  
Peach introduced rhombal algebras in his thesis~\cite{P} by  
imposing certain relations on quivers corresponding to  
tilings of the plane by rhombi. He shows that quotients of these 
rhombal algebras model parts of weight 2 blocks of symmetric groups. 
Ringel~\cite{Ri} and Turner~\cite{T} have further  
analysed these rhombal algebras, and   
Chuang and Turner~\cite{CT} generalised them to higher dimensions.  
Among the general methods to analyse infinite-dimensional  
rhombal algebras we develop, the most important tool  is a filtration by rhombus modules. In particular, 
to each rhombus in the quiver we associate a module, which we call a rhombus module.  We then determine  
how the projective indecomposable modules and others  
are filtered. In addition to more precise 
information on the module structure of the modules analysed, these rhombus filtrations provide  
insight into the structure of the module category.

We apply our methods to the infinite-dimensional rhombal algebra  
associated to a special rhombal tiling known as the Rauzy tiling.  
The Rauzy tiling is constructed in stages by a substitution rule that reflects the 
underlying self-similarity of the tiling. The dynamics and combinatorics 
of the Rauzy tiling have been extensively studied (see, for example,~\cite{ABS},\cite{BV} or \cite{F}), 
and we make use of its special properties for  
the solution of many of the representation theoretic problems we encounter. 
 
More precisely, we develop  techniques relating the 
infinite-dimensional Rauzy algebra to  finite-dimensional  
algebras related to the finite portions of  
the complete tiling that occur in its construction.  
Comparing different possible truncation  
methods for making such relations leads to a natural choice 
of truncating the infinite-dimensional algebra by an idempotent  
associated to the finite portions, resulting in what we call the finite-dimensional Rauzy algebras. 
Applying our  technique of rhombus filtrations to analyse Rauzy algebras  
leads to an understanding of many of the modules that arise in this setting. 
Special features of the Rauzy tiling and our general techniques allow us to reveal  
many interesting properties of the finite-dimensional Rauzy  
algebras. 
 In particular, we are able to  
identify modules with periodic orbits under application of the  
Heller operator. These periodic orbits occur along lines in the quiver that 
are analogues of the `induced lines' introduced in~\cite{EM}. Modules 
with periodic orbits are crucial for understanding the graph 
structure of the stable Auslander-Reiten quivers of the finite-dimensional 
Rauzy algebras. Through our results we are able to identify 
the tree classes and hence the graph structure of many components. Some of the features which we  
found lead to a general conjecture, 
including the occurrence of trees of class  $A_{\infty}^{\infty}$.

\section{The Rauzy tiling} \label{tiling}  
Rauzy~\cite{R} introduced a fractal domain arising in a natural  
way from a substitution. Both the fractal  
and the substitution now bear his name.  
To better understand the Rauzy fractal and its dynamics,  
alternative geometric interpretations have been introduced. As  
shown in~\cite{ABS},~\cite[Chapt. 8]{F}, there is an  
increasing sequence of patches $\p_i$ of tiles whose union  
$\bigcup \p_i$ forms a complete tiling of the plane. Appropriate  
renormalizations of the patches $\p_i$ converge to the Rauzy  
fractal, but we shall only be interested in the tiling of the  
plane, which we refer to as the \textit{Rauzy tiling}. We shall  
describe a recursive method for obtaining the increasing sequence  
of patches $\p_i$, the edges and vertices of which yield the  
quivers used to construct the Rauzy algebra. This recursive  
construction reflects the original substitution and produces a  
tiling with a hierarchical, self-similar structure.  
  
The tiles in $\p_i$ are  projections of faces of unit  
cubes in $\mathbb{R}^3$ with vertices in $\mathbb{Z}^3$.  
Specifically, we define the linear map  
$\mathrm{p}:\mathbb{R}^3 \rightarrow \mathbb{R}^2$ by  
\[  
\mathrm{p}\left(\mathbf{e}_1\right)=\left(-\frac{\sqrt{3}}{2},-\frac{1}{2}\right),\:\mathrm{p}  
\left(\mathbf{e}_2\right)=\left(\frac{\sqrt{3}}{2},-\frac{1}{2}\right),\:  
\text{and } \mathrm{p}\left(\mathbf{e}_3\right)=\left(0,1  
\right)=-\left(\mathrm{p}\left(\mathbf{e}_1  
\right)+\mathrm{p}\left(\mathbf{e}_2\right)\right),  
\]  
where the $\mathbf{e}_i$ are the standard basis  
elements of $\mathbb{R}^3$. Then $\mathrm{p}$ projects the lattice $\mathbb{Z}^3$ onto the  
planar lattice  
\[  
\mathcal{L}=\left\{ m\, \mathrm{p} \left(\mathbf{e}_2\right)+n\,  
\mathrm{p}\left(\mathbf{e}_3\right) \: : \: m,n\in  
\mathbb{Z}\right\},  
\]  
and all tiles will be rhombi with vertices in $\mathcal{L}$.

\begin{figure}[ht]  
\centering \epsfig{file=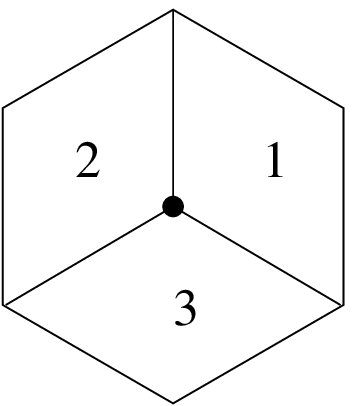,  height=2cm} \caption{$\p_0$}  
\label{step0}  
\end{figure}  
  
In fact, each rhombus of the tiling is a translation of one of the three  
rhombi $R_1,\,R_2$ or $R_3$  of the initial patch $\p_0$  
pictured below, where the common point of all three is the origin.  
We shall describe a rhombus in $\p_i$ as a  
translation of one of the original rhombi of the form $\mathrm{p}  
\left(z\right)\,+\,R_k$ for some $z \in \mathbb{Z}^3$. As $\mathrm{p}$ is not injective,  
this representation is not unique. However, to each rhombus we will 
recursively assign a specific representation of this form, where we identify   
$R_k$ with $\mathrm{p} \left({\mathbf 0} \right)\,+\,R_k$. The  
recursive process in going from $\p_i$ to $\p_{i+1}$ is determined  
by the function $\mathcal{R}$ from the set of rhombi in $\p_i$ to the 
collection of rhombi in $\p_{i+1}$. 
A key element in the definition of $\mathcal{R}$ is the matrix $\mathbf{M}=\left(%
\begin{array}{rrr}  
  0 & 1 & 0 \\  
  0 & 0 & 1 \\  
  1 & -1 & -1 \\  
\end{array}%
\right)$  derived from the Rauzy substitution.  We denote the  
columns of $\mathbf{M}$ by $\mathbf{c}_i$, considered as elements  
of $\mathbb{Z}^3.$ Then one obtains the rhombi in $\p_1$ by taking the union of the following sets of rhombi:  
\[  
\mathcal{R}\left(R_1\right)= \left\{R_3 ,\, \mathrm{p}  
\left(\mathbf{c}_2 \right)+ R_1 ,\, \mathrm{p}  
\left(\mathbf{c}_3\right)+\,R_2  
\right\},\:\mathcal{R}\left(R_2\right)= \left\{ R_1\right\},\:\text{and  
}\mathcal{R}\left(R_3\right)= \left\{ R_2\right\}  
\]  
In general, once the rhombi in $\p_i$ have been determined, one defines on each rhombus  
$\mathrm{p} \left(z \right) +\,R_k \in \p_i$  
\[ 
\mathcal{R}\left(\mathrm{p} \left(z \right)  
+\,R_k\right):=\mathrm{p}  
\left(\mathbf{M}z\right)+\,\mathcal{R}\left(R_k\right),  
\]  
where the translation by $\mathrm{p} \left(\mathbf{M}z\right)$  
applies to each rhombus in $\mathcal{R}\left(R_k\right)$, and one represents $\mathrm{p}  
\left(\mathbf{M}z\right)+\,\mathrm{p}  
\left(\mathbf{c}_{k+1}\right)+\,R_k$ as $\mathrm{p}  
\left(\mathbf{M}z + \mathbf{c}_{k+1}\right)+\,R_k$. One  
then obtains the rhombi in $\p_{i+1}$ by taking the union of all the rhombi in  
$\mathcal{R}$ applied to  the rhombi in $\p_{i}$.  
  
\begin{figure}[ht]  
\centering  
\epsfig{file=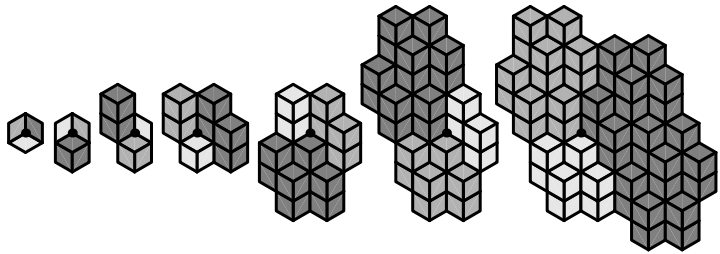, width=10cm}  
\caption{$\p_0$ through $\p_6$, with  $\mathcal{R}^i\left(R_k\right)$ in $\p_i$ in darkness descending as $k$ increases.}  
\label{Steps0-4}  
\end{figure}

  

\section{Rhombal algebras associated to the Rauzy tiling}

\subsection{Overview of path algebras of quivers} 
 
We survey the  properties of the representation theory of quivers and   
paths algebras of quivers  vital to our investigations. For a more detailed reference,  
we refer the reader to \cite{ARS} and \cite{B}.  
 
\medskip 
 
A quiver $Q$ is an oriented graph. A path in the quiver is a sequence of arrows such that each arrow  
begins at the 
vertex at which the preceding arrow ended. Furthermore, to each vertex $z$ we associate a trivial path $e_z$.  
For a field $k$, the path algebra $kQ$ of the quiver  is the algebra generated by the set of all paths of  
$Q$, where the multiplication of two paths is given by concatenation of the two paths if they concatenate and 
by zero otherwise. Note that for each vertex $z$ of $Q$, $e_z$ is an idempotent in $kQ$. 
 
Given a set of relations $R$ on the paths of $Q$, we can define a two sided ideal $I$ generated by $R$,  
leading to the quotient algebra $kQ/I$.  
Then by a fundamental theorem of Gabriel~\cite{G}, every basic algebra over an algebraically closed field  
is isomorphic to the  path algebra of a quiver with relations.   
 
The radical of $kQ/I$ is spanned by the image of all paths of length $\geq 1$. The socle of $kQ/I$ 
is spanned by all paths $b$ such that $ab=0$ for all arrows $a$. 
 
Furthermore, for every vertex $z$ of $Q$, there is a simple $kQ/I$-module $S_z$ of dimension one associated to $z$  
such that $S_ze_z = S_z$ and $ S_ze_x =0$ for all $x \neq z$, and for any arrow of $Q$ the corresponding element 
in $kQ/I$ acts as zero. Then $e_z(kQ/I)$ is the projective indecomposable  
corresponding to $S_x$. For each vertex $x$, the dimension of $e_z(kQ/I)e_x$ is equal 
to the composition multiplicity of $S_x$ as a composition factor of $e_z(kQ/I)$.

\subsection{Infinite-dimensional Rauzy algebra}\label{infinite}  
\begin{figure} [ht] 
\centering 
\input{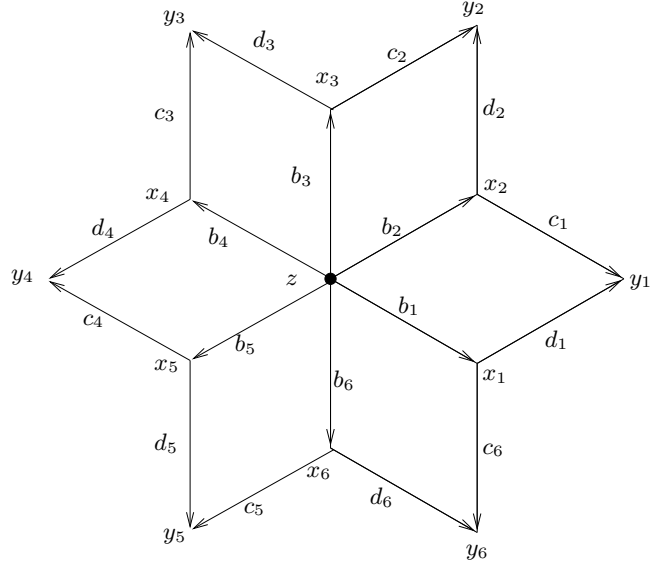} 
\caption{Labelling of the vertices and paths neighbouring $z$} 
\label{star6} 
\end{figure}

The graph formed by the edges and vertices of the Rauzy tiling is a member of the family of  
graphs investigated by Peach~\cite{P}.   
We form a quiver $Q$ from the Rauzy  tiling  by  
replacing each edge with a double arrow in opposite directions. In representing quivers, letters towards the end of the alphabet are reserved for  vertices, while letters towards the beginning 
of the alphabet are used for arrows. When focusing on a particular vertex $z$, arrows in figures will be drawn with tail  
given by the vertex that is closest to $z$, and arrows and vertices around $z$ will be labelled with indices increasing  
in the counter-clockwise direction, as indicated in figure~\ref{star6}. The arrow pointing in the opposite direction  
from the arrow $b$ is denoted $\bar{b}$.

\bigskip 
 
For the path algebra $kQ$ associated to the field $k$,   
following Peach~\cite{P} we define a set of relations $R$ on $kQ$ as follows:  
  
1. (Two rhombus relation) Any path of length two in $kQ$ that borders more than one rhombus is zero.  
  
2. (Mirror relation) Any two paths of length two connecting opposite vertices of the same rhombus are equal.  
  
3. (Star relation) For any vertex $z$ and labelling as indicated in figure~\ref{star6}, the following relations hold, where we replace a path by zero if there is no  corresponding edge in the quiver.  
$$  
\begin{array}{ccc}  
b_2\bar{b_2}-b_5\bar{b_5} &=& \overset{|}\varepsilon_z (b_4\bar{b_4}-b_1\bar{b_1})\\  
b_6\bar{b_6}-b_3\bar{b_3} &=&  \overset{\diagdown}\varepsilon_z (b_2\bar{b_2}-b_5\bar{b_5})\\  
b_4\bar{b_4}-b_1\bar{b_1} &=& \overset{\diagup} \varepsilon_z (b_6\bar{b_6}-b_3\bar{b_3})\\  
\end{array}$$ 
where $  \overset{|}\varepsilon_z, \overset{\diagdown}\varepsilon_z, \overset{\diagup}\varepsilon_z \in \{ -1,1\} $ and  
$ \overset{|}\varepsilon_z \overset{\diagdown}\varepsilon_z \overset{\diagup}\varepsilon_z =1$. 
A choice of signs at a single vertex $z$ determines the signs at all vertices of $Q$. For more detail we refer the reader to~\cite{P}. 
 
\bigskip

\begin{definition} {\rm The  infinite-dimensional} {\it Rauzy algebra} $A$  {\rm is defined to be}  $kQ/I$, {\rm where} 
 $I$ {\rm is the  
ideal generated by the relations} $R$.  
 \end{definition} 
  
\bigskip 
 
Let $z$ be a vertex and $P_z = e_z A$ the corresponding projective module of $A$. Peach ~\cite[Ch. 2]{P} has 
shown many facts about paths in $A$ and  $P_z$ that we list here for convenience.  
 
\medskip 
 
Any path which does not border a single rhombus is zero in $A$. Furthermore, 
all paths of length $\geq 5$ are zero in $A$, and all paths of length 4 
which do not start and end at the same vertex are zero in $A$.  
The socle of $P_z$ is simple and isomorphic to $S_z$. In fact, $A$ is a symmetric algebra. 
The image of any path $p$ of length four around a rhombus starting and ending at 
$z$ spans the socle of $P_z$.  
Any two paths of length 3 around a single rhombus are linearly dependent. 
In particular, if $x$ is a vertex such that there is an arrow $b$   
from $z$ to $x$, 
then the space $e_zAe_x$ has dimension $2$, and it is spanned by $b$ together 
with any path of length three around a rhombus which has corners $z, x$.

\subsection{Finite-dimensional Rauzy algebras}\label{finalg}

Each patch $\p_i$ in the formation of the Rauzy tiling as described in  
section \ref{tiling} gives rise to a quiver $Q_i$  
with vertex set $V_i$\,. Since $\p_i \subset \p_{i+1}$ for each $i$ and since the Rauzy tiling itself is given by  
$\cup \p_i$, the underlying graph of each $Q_i$ embeds naturally in the infinite graph given by the Rauzy  
tiling. Two natural possibilities then arise for relating  the finite-dimensional path algebras  associated  
to the $Q_i$  with the Rauzy algebra $A$. We can either truncate the vertices but allow all paths that start   
and end in a given $V_i$, or we can truncate the relations together with the vertices. 
 
\medskip 
 
\begin{figure}[ht]  
\centering  
\input{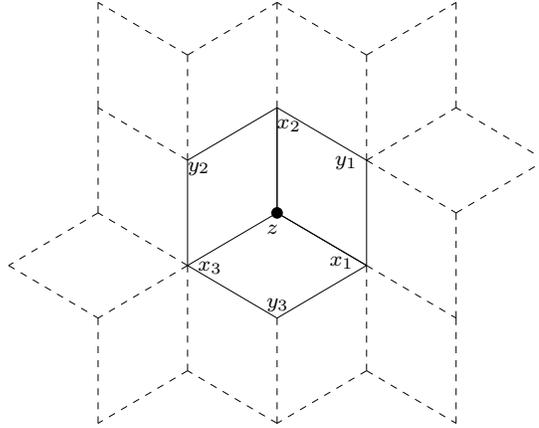} 
\caption{The vertices $x_i$, $y_i$ and $z$ are the vertices of $Q_0$}  
\label{Step0}  
\end{figure}

The first option corresponds to  
cutting $A$ by idempotents given by the sum of the primitive idempotents of the vertices in $V_i \,$.  
With this approach, $Q_i$ yields the algebra 
\[   
A_i = e_i A e_i, \mbox{where $e_i = \sum_{v \in V_i} e_v$ and $e_v$ is the   
primitive idempotent at $v$}. 
\] 
Thus, there will be paths which are non-zero  in $A_i$ that go through a vertex not contained in $V_i \,.$  
The dotted lines in figure~\ref{Step0} indicate possible edges for non-zero paths of $A_0$  
starting and ending at vertices in $V_0$   
which are not paths of  $Q_0\,$.  
  
\bigskip  
  
\medskip 
 
The second option corresponds to a  
truncation of the  underlying graph on which we then impose the relations $R$, yielding 
the algebras $B_i = kQ_i /(I \cap kQ_i)$  with the relations $R$. 
  
To decide which of these  approaches to choose in order to define  
the finite-dimensional path algebras, we examine the indecomposable projective modules.  
  
\subsection{Which type of truncation?}

Let $z$ be a vertex and $P_z = e_zA$ the corresponding indecomposable projective module of $A$,  
where $e_z$ is the primitive idempotent at $z$.  
The star relations imply that the crucial portion of $P_z$ is the space $e_zP_ze_z$. 
More precisely, the following subspace of $e_zP_ze_z$ and $P_z$ proves to be most important.   
 
\begin{definition} 
{\rm Let $X_z$ be the subspace of $P_z$ spanned by paths of length two starting and ending at $z$.} 
\end{definition}

\bigskip

We calculate the dimension of $X_z$ as a subspace of $A$ and as a subspace of the finite-dimensional  
algebras $A_i =e_iAe_i$. In fact, we show that in the latter case, $X_z$ is also a subspace of a  
projective module at $z$ if the primitive idempotent at $z$ is a summand of the idempotent $e_i$.    
We wil also show that imposing the star relations on $B_i$ 
produces a finite-dimensional algebra which is not symmetric,  which we will therefore not  
consider in later sections.   
  
\bigskip  
  
The key difference in  the relations for the two types of truncations lies in the star  
relations.  
  
\bigskip  
  
{\it (I) Cutting $A$ by the idempotent $e_i$.} Recall that $e_i$ is the sum of the primitive  
idempotents in $Q_i$ and that $A_i = e_i Ae_i$.   
In this case the star relations are as given in section~\ref{infinite}.  
Let $z$ be a vertex   
such that $e_ze_i=e_z$. Then $e_z$ is  
a primitive idempotent of the algebra $e_iAe_i$, and we have the   
indecomposable projective modules $P_z = e_zA$ of $A$, and  
$e_zAe_i = e_z(e_iAe_i)$ of $e_iAe_i$.   
We now compare $X_z$ with the subspace of $e_zAe_i$ given by paths of length two.    
  
 
Since $e_ie_z=e_ze_i=e_z$, we have $e_zAe_ze_i = e_zAe_z$, and as a vector space $e_zAe_z$ decomposes into  
$e_zAe_z = \langle e_z\rangle \oplus X_z \oplus \langle X_z^2 \rangle$. It follows from  
\cite[2.4.12]{P} that the last summand   
is $1$-dimensional. The decomposition still holds when we multiply by $e_i$ on both sides.  
Therefore, in $e_zAe_i$ the space spanned by paths of length two beginning and ending at $z$  
is equal to $X_z$,  and $X_z$ has the   
following dimensions:    
  
$$\dim(X_z) = \left\{\begin{array}{ll} 1 & \mbox{$z$ is a 3-vertex}  
\cr  
2  & \mbox{$z$ is a 4-vertex }\cr  
3  & \mbox{$z$ is a 5-vertex }\cr  
4  & \mbox{$z$ is a 6-vertex }  
\end{array}  
\right.  
$$  
 
\bigskip  
  
To examine the second type of truncation, the following definition is essential.  
\begin{definition} 
{\rm Let $z$ be a vertex of  $Q_i$. We say that $z$ is a {\it $(k,n)$-vertex} 
if $z$ is an $n$-vertex in the untruncated quiver $Q$ and a $k$-vertex in  $Q_i$.}  
\end{definition}

\begin{figure}[ht]  
\centering  
\input{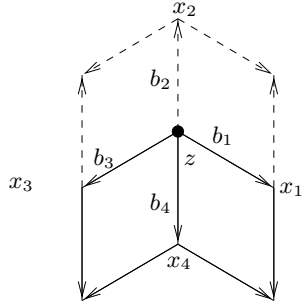}  
\caption{(4,3)-vertex}  
\label{vertex3-4}  
\end{figure}

For example, $z$ in figure~\ref{vertex3-4} is a $(3,4)$-vertex, where the dotted arrows are outside $Q_i$. 
  
We label  the patch in the quiver around $z$ as indicated in figure~\ref{star}, with the same condition 
on arrows as before. 
\begin{figure} [ht] 
\centering 
\input{starA.pstex_t} 
\caption{} 
\label{star} 
\end{figure} 
Then if for example $z$ is a three vertex, the patch in the quiver around  
$z$ will be labelled as in figure~\ref{threevertex}. 
\begin{figure} [ht] 
\centering 
\input{threevertexA.pstex_t} 
\caption{} 
\label{threevertex} 
\end{figure} 
  
\bigskip  
  
{\it (II) Imposing the relations on $Q_i$.}  
Recall that $B_i = kQ_i /(I \cap kQ_i)$; that is,  we impose the 
relations $R$ on the truncated quiver. If $z$ is a vertex of $Q_i$, we calculate  
the dimension of $X_z$ as a subspace of $B_i$.    
  
(a) Suppose $z$ is a $(2,5)$-vertex in $Q_i$.   
Then $X_z$ is spanned in $A$ 
by two paths, and one pair of opposite paths is completely missing.  
Also, the opposite of each existing path is missing. If we now impose the star  
relation on the truncated quiver $Q_i$, all of the two paths in $X_z$ are actually equal to zero.

\medskip  
  
(b) Suppose $z$ is a $(6,3)$-vertex in $Q_i$.  
Then as a subspace of $A$, $X_z$  is spanned by three paths.  
If we impose the star relations on $Q_i$, 
then any two of these paths are linearly dependent, and thus $X_z$ as a subspace of $B_i$  is 1-dimensional.

\medskip  
  
(c) Suppose $z$ is a $(4,3)$-vertex in $Q_i$. 
 
With the notations of figure~\ref{vertex3-4} the star relations in $A$ for such a vertex are   
$$\begin{array}{ccc}  
-b_3\bar{b}_3 &=& {\overset{|}\varepsilon}_z(-b_1\bar{b}_1) \\   
b_4\bar{b}_4-b_2\bar{b}_2 &=& \overset{\diagdown}\varepsilon_z(-b_3\bar{b}_3) \\    
-b_1\bar{b}_1 &=& \overset{\diagup}\varepsilon_z(b_4\bar{b}_4-b_2\bar{b}_2)\\  
\end{array}  
$$  
However, in the truncated quiver $b_2$ does not exist. So the imposed star relations give that any two of the   
three paths are linearly dependent and $X_z$ as a subspace of $B_i$ is again $1$-dimensional.   
 
\medskip 
 
As a special case, we consider the following.  
  
\medskip

{\it Imposing star relations on $Q_1$.}

\bigskip  
  
\begin{figure}[h!]  
\centering  
\input{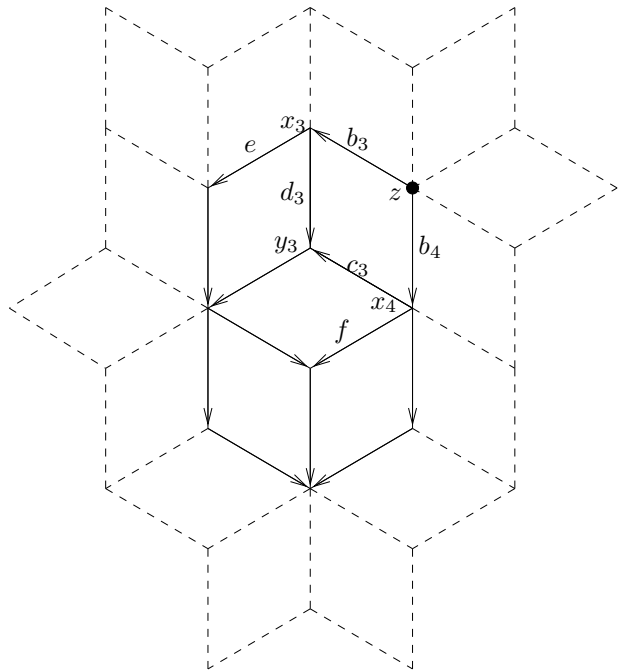}  
\caption{The solid arrows are in $B_1$ and the dotted lines represent  
arrows that are in $A$ but not in $B_1$. The arrows are 
defined with respect to $z$, and we have only indicated at most one direction for each arrow.} 
\label{Step1}  
\end{figure}

Let  $z$ be the  $(5,2)$-vertex of $Q_1$. Then with the notation established in 
section 3.2, we get a labelling as in figure~\ref{Step1}. 
 
\bigskip  
 
Now consider the space $P_z':=e_z(B_1)$.   
We have seen that for the vertex $z$, the paths $b_3\bar{b}_3=0=b_4\bar{b}_4$,  
where $b_3, \bar{b}_3, b_4$ and $\bar{b}_4$ are   
paths inside $B_1$. By the mirror relation,  
$b_3d_3 = b_4c_3$. Thus,   
for a basis of $P_z$ the only path of length two we need is  $b_3d_3.$  
Now consider paths of length three. There are only two such paths to consider: $b_3d_3\bar{d}_3$ and $b_3d_3\bar{c}_3$  
(all others are zero).  
But $d_3\bar{d}_3$ occurs in a star relation starting at $x_3$, and from II(c) above we know  
that any two paths of length two are dependent. So $d_3\bar{d}_3 = e\bar{e}$, where $e \neq d_3$ is some path  
starting at $x_3$. But then $b_3d_3\bar{d}_3 = b_3e\bar{e} =0$.     
  
By the mirror relation  $b_3d_3\bar{c}_3 = b_4c_3\bar{c}_3$. We have an imposed star relation  
at $x_4$ of the form  
$$(0-f\bar{f}) = \overset{|}\varepsilon_{x_4}(c_3\bar{c}_3 - 0).  
$$  
Hence, by the two rhombus relation, $b_4c_3\bar{c}_3 = \pm b_4f\bar{f} = 0$.   
  
\bigskip  
  
This has proved that $(b_3d_3)J=0$ (with $J= \mathrm{rad}(e_1Ae_1)$, and thus  
$b_3d_3$ lies in the (right) socle. So the projective $P_z'$  has  Loewy length  
at most $3$ (most likely, it is $3$); and, moreover, the socle is not isomorphic to the top.  
Therefore, this is not a projective module for a symmetric algebra!  
(Although, the algebra is quite possibly self-injective).

\bigskip

{\bf Conclusion~:} Imposing the star relations on the truncated quiver does not produce  
symmetric algebras, and hereafter the only finite-dimensional algebras we consider are  
the $A_i$.

\section{Rhombus modules and rhombus filtrations}

In this section we develop a general technique that applies 
to all infinite-dimensional rhombal algebras. We describe how to associate modules, which we call 
rhombus modules,  to  any rhombus of the quiver associated to the algebra.     
These modules may be thought of as building blocks of the rhombal algebras.  
We show that the projective modules of the algebra $A$   
have filtrations by such rhombus modules.  
   
\bigskip  
   
\subsection{Rhombus modules}  
  
Consider a rhombus in the quiver labelled as in figure~\ref{rhombus1}. 
  
\begin{figure} [ht] 
\centering 
\input{rhombus1.pstex_t} 
\caption{} 
\label{rhombus1} 
\end{figure} 
  
\bigskip

\vspace*{0.2cm}  
  
\begin{lemma}\label{rhombusmodule}  For any rhombus in the quiver,  
there is a unique module $R$ with top isomorphic to $S_y$ and socle  
isomorphic to $S_z$ and $rad(R)/soc(R)\cong S_{x_1}\oplus S_{x_2}$.  
Explicitly, $R=b_1dA = b_2cA$. 
\end{lemma}

\bigskip

\begin{definition} 
{\rm We say that a module $R$ such as in Lemma~\ref{rhombusmodule} is a {\it rhombus module}  
and denote it by $R=R^y_z$.}   
\end{definition} 
  
\bigskip

{\it Proof } (1) Existence: The space  $e_zAe_y$ is the 1-dimensional space  
spanned by $b_1d$. The module $R:= (b_1d)A$ is contained in $e_zA$ and therefore  
has socle $S_z$. The top of $R$ is isomorphic to $S_y$ as it  
is generated by an element $x$ with $x=xe_y$. Furthermore, $R$ has basis  
$$b_1d,   b_1d\bar{d},  b_1d\bar{c}, b_1d\bar{d}\,\bar{b}_1  
$$  
That the elements listed are all non-zero follows from ~\cite[2.4.13]{P}.

(2) Uniqueness: Any such module must be isomorphic to a submodule  
of the injective module with socle $S_z$; that is, a submodule of $e_zA$.  
As it is generated by an element of $e_zAe_y$ and this is 1-dimensional and  
spanned by $b_1d$, it follows that $R\cong b_1dA$.  \hfill $\Box$

\bigskip  
  
{\bf Remark }   
There are four such rhombus modules for each rhombus in the quiver. (This may be considered as an analogue of  
the modelling of different  $\mathcal{O}$-forms for a liftable module with an irreducible character.) 
 
\bigskip

\subsection{Rhombus filtrations}  
  
We shall show that several modules that arise naturally in our setting,  
including the modules generated by arrows and indecomposable modules,   
have rhombus filtrations as defined below. 
  
\begin{definition} 
{\rm We say that a module $M$ has a {\it rhombus filtration} if there is a sequence  
of submodules $0=M_0 \subset M_1\subset \ldots M_k = M$,  
where $M_i/M_{i-1}$ is a rhombus module for each $i$.} 
\end{definition}   
 However, as the following lemma shows,  
`rhombus filtration multiplicities' are not well-defined if one allows 
all rhombus modules as quotients.  
 
\medskip 
  
Throughout this section we refer to figure~\ref{star} for notation.  
  
\medskip  
  
\begin{lemma} {\rm (Arrow Lemma)} Suppose $b$ is an arrow in the quiver. Then   
the module $bA$ has two  rhombus filtrations, each with   
two quotients.  
\end{lemma} 
{\it Proof } Use the standard labelling of vertices and arrows as defined in  
figure~\ref{star} and set $b=b_1$.   
 
(1) There is an exact sequence   
 
$$0\to bc_nA\to bA \to \bar{b}_2bA\to 0.  
$$ 
Clearly $bA$ contains the rhombus module $bc_nA$, and furthermore left   
multiplication by $\bar{b}_2$ induces a surjection $\pi: bA \to    
\bar{b}_2bA$. We also have $\bar{b}_2bc_n=0$, so   
$bc_nA \subseteq {\rm ker \,}\pi$. Then equality holds by dimensions, since  
$bA/(bc_n)A$ has basis the cosets of  
$b, bd_1, b\bar{b}, bd_1\bar{c}_1$.  
  
(2) Similarly there is an exact sequence  
 
$$0\to bd_1A\to bA \to\bar{b}_nbA\to 0.  
$$ 
\hfill $\Box$

Next we consider submodules of projectives which are generated  
by two arrows.

\begin{lemma} {\rm (Two-arrow lemma)}  With the notation of figure~\ref{star}, assume  
$z$ is an $n$-vertex and $n\geq 4$.   
Then the module $b_1A + b_2A$ has a rhombus filtration, with three rhombus  
quotients. Explicitly, we have  
an exact sequence  
 
$$0\to b_1d_1A \  \to \   
b_1A + b_2A \to  \bar{b}_nb_1A \oplus \bar{b}_3b_2A \to 0  
 $$  
\end{lemma} 
{\it Proof }  Since $b_1d_1=b_2c_1$, we have a commutative diagram with exact rows and columns  
  
 $$\CD 
&& 0 && 0 && 0 &&\\  
 &&@VVV @VVV @VVV && \\ 
0@>>>{\rm ker \,} p_1 @>>>{\rm ker \,} p@>>>{\rm ker \,}\phi@>>> 0\\  
 &&@VVV @VVV @VVV && \\ 
 0@>>> b_1d_1A\oplus b_2c_1A @>>> b_1A\oplus b_2A @>>> \bar{b}_nb_1A\oplus \bar{b}_3b_2A@>>> 0\\ 
 &&@Vp_1VV@VVpV@VV\phi V \\ 
 0@>>>b_1d_1A@>j>>b_1A+b_2A @>>>C@>>>0\\ 
 &&@VVV @VVV @VVV && \\ 
&& 0 && 0 && 0 &&\\ 
 \endCD 
 $$

 Namely, take as the middle row the direct sum of the two arrow sequences.  
 Then take for $p$ and $p_1$ the addition maps, and the map $j$ is inclusion.  
 Then the left lower square commutes, and hence it induces the map $\phi$ making  
 the right lower square commute as well. Clearly, $p$ and $p_1$ are surjective,  
 and then $\phi$ is also surjective. Thus, the top row is exact by  
the  Snake Lemma.  
   
By definition, we have ${\rm ker}\, p = \{ (x, -x): x \in b_1A\cap b_2A \}$,  
which is isomorphic to $b_1A\cap b_2A$. Similarly,  
${\rm ker}\, p_1 = (b_1d_1, -b_2c_1)A \cong b_1d_1A$, which is a rhombus module.

\medskip  
 
If $n\geq 4$, we now show that $b_1A\cap b_2A = b_1d_1A$. 
To see this, note that $b_1d_1 = b_2c_1$, and  $b_1d_1A$ is thus contained in the intersection. 
 
Suppose to the contrary that the intersection were not contained in $b_1d_1A$. Then there would be 
a simple module contained in both $b_1A/b_1d_1A$ and $b_2A/b_2c_1A$. 
From the Two-arrow Lemma, we know that 
$b_1A/b_1d_1A \cong \bar{b}_nb_1A$ and that it has simple socle 
isomorphic to $S_{x_n}$. Similarly, $b_2A/b_2c_1A\cong \bar{b}_3b_2A$, which 
has simple socle isomorphic to $S_{x_3}$.  
For $n\geq 4$ the vertices $x_n$ and $x_3$ are distinct, and then 
$S_{x_n}$ is not isomorphic to $S_{x_3}$, a contradiction that demonstrates 
the desired equality.

 Hence in this case, ${\rm ker}\, \phi = 0$ and $\phi$ is therefore an isomorphism.   
\hfill $\Box$ 
\medskip  
  
{\bf Remark. } The condition that $n\geq 4$ is necessary.   
If $z$ is  
a 3-vertex, then the intersection of $b_1A$ and $b_2A$ properly contains  
the rhombus $b_1d_1A$ and $\phi$ has a kernel of dimension two.

\medskip

We now prove a `dual' version of the previous lemma. As the arrow $\bar{b}_1$ starts at the vertex 
$x_1$ and the arrow $\bar{b}_2$ starts at the vertex $x_2$, $(\bar{b}_1, \bar{b}_2)A$ is a submodule of 
$e_{x_1}A\oplus e_{x_2}A$. 
 
\medskip

\begin{lemma} {\rm (Dual two-arrow lemma)}\label{Dual two-arrow lemma}    
 Suppose $z$ is an $n$-vertex with $n\geq 4$. Then   
$(\bar{b}_1, \bar{b}_2)A$ has a rhombus filtration with three rhombus  
quotients. Explicitly, there is an exact sequence  
$$0 \to \bar{b}_1b_nA \oplus \bar{b}_2b_3A \to   
(\bar{b}_1, \bar{b}_2)A \to \bar{d}_1\bar{b}_1A\to 0  
$$  
\end{lemma}

\bigskip  
 {\it Proof }  We have a commutative diagram with exact rows and columns  
 
$$\CD  
&& 0 && 0 && 0 &&\\  
 &&@VVV @VVV @VVV && \\ 
 0@>>>{\rm ker \,} \pi_1 @>>>(\bar{b}_1, \bar{b}_2)A@>\pi_1>>(\bar{d}_1\bar{b}_1, \bar{c}_1\bar{b}_2)A@>>> 0\\ 
  &&@V\phi VV @Vj VV @Vj_1 VV  \\ 
   0@>>> \bar{b}_1b_nA\oplus \bar{b}_2b_3A @>>>  
   \bar{b}_1A\oplus \bar{b}_2A @>\pi>> \bar{d}_1\bar{b}_1A\oplus \bar{c}_1\bar{b}_2A@>>> 0\\ 
    &&@VVV@VVV@VVV \\ 
     0@>>>{\rm Cok\,}\phi @>>>{\rm Cok\,}j@>>>{\rm Cok\,}j_1@>>>0\\ 
 &&@VVV @VVV @VVV && \\ 
&& 0 && 0 && 0 &&\\ 
      \endCD 
       $$ 
 
Namely, start with the middle row which is a direct sum of the two arrow sequences.  
Then take for $j$ and $j_1$ the inclusion and for $\pi_1$ the restriction of $\pi$.  
Then the top right square commutes and  
induces the map $\phi$, which then also is 1-1. Now the Snake Lemma gives  
the lower row. Since $\bar{d}_1\bar{b}_1A = \bar{c}_1\bar{b}_2A$, it follows that  
${\rm Coker \,} j\cong \bar{d}_1\bar{b}_1A$ which is a rhombus module.   
If $n\geq 4$, then we check that ${\rm Cok\,} j \cong {\rm Cok\,}j_1$  
(for example, by dimensions), and then $\phi$ is an isomorphism.  
\hfill $\Box$ 
  
As before, when $n=3$ this does not hold.

\bigskip

\begin{lemma} {\rm (Three-arrow lemma)}  Suppose $n\geq 5$. Then the module  
$M=b_1A+b_2A+b_3A$ has a rhombus filtration with four quotients. Explicitly, we have an exact sequence  
$$0\to b_1A+b_2A \to M \to \bar{b}_4b_3A\to 0.  
$$ 
\end{lemma}

\bigskip  
  
{\it Proof } The module $M$ is contained in $e_zA$. Define   
$\pi: M\to \bar{b}_4A$ to be multiplication by $\bar{b}_4$, a homomorphism  
from $M$ onto $\bar{b}_4b_3A$. From~\cite{P} (see the  
summary in section~\ref{infinite}), we know that  
$b_1A+b_2A$ is contained in the kernel of $\pi$. It remains to show that  
equality holds.   
  
Restrict $\pi$ to $b_3A$, the Arrow Lemma shows that  
${\rm ker}\, \pi \cap b_3A = b_3c_2A$. But $b_3c_2 = b_2d_2 \in b_2A$, and so the intersection  
of ${\rm ker}\, \pi$ with $b_3A$ is contained in $b_1A+b_2A$. Now it follows that  
${\rm ker}\, \pi \subseteq b_1A + b_2A$.   
The claim follows. \hfill $\Box$

Below is a dual version. The proof, which is similar, is left to the reader.  
\bigskip  
  
\begin{lemma}{\rm (Dual three-arrow lemma)}  Let $n\geq 5$. The module $M:= (\bar{b}_1, \bar{b}_2, \bar{b}_3)A$ has a rhombus filtration with four quotients, and we have  
an exact sequence  
$$0\to \bar{b}_3b_4A \to M \to (\bar{b}_1, \bar{b}_2)A\to 0. 
$$ 
\end{lemma}

\subsection{Indecomposable projective modules}

\bigskip

First we  observe that the composition factors of projective modules have an easy description in  
terms of   possible rhombus filtrations.   For this  we use the following standard notation from  
the representation theory of quivers: for an  $A$ module 
$M$, $\underline{\dim}\; M = (m_x)_{x}$, where $x$ ranges over the set of vertices of $Q$ 
and $m_x = \dim\, Me_x$ for  the primitive idempotent $e_x$ at $x$.  
 
Recall that $[e_zA:S_x] = \dim e_zAe_x$. For this to be nonzero, 
there must be a rhombus which has corners 
$z$ and $x$. As noted in  \ref{infinite}, if there is an arrow from $z$ to $x$, the dimension is $2$.  
If $z$ and $x$ are opposite corners, then the dimension is 
$1$; specifically, any two paths of length 2 are the same by the mirror relation, and 
there are no non-zero paths of length $> 2$.  
Furthermore, the dimension of $e_zAe_z$ is $n$ if $z$ is an $n$-vertex.  
Together these facts imply the following lemma. 
  
\begin{lemma}{\rm  (Multiplicity Lemma)}  The   
composition factors of $e_zA$ are given by   
$$\underline{\dim}\, e_zA = \sum_i \underline{\dim}\, (R^z_{a_i}).  
$$  
That is, $[e_zA:S_z] = n$ where $z$ is an $n$-vertex,   
$[e_zA:S_{x_i}] = 2$ and $[e_zA:S_{y_i}]=1$, for each $i$.  
\end{lemma} 
  
The same applies to the modules $e_zA_i$ in $A_i$ for $z$ sufficiently far from the boundary of $Q_i$.  
  
\bigskip  
  
We now describe rhombus filtrations of indecomposable projectives with the aid of the 
the {\it Heller operator} $\Omega$. Recall that for a module $M$ with projective cover $p: P_M \to M$, 
 $\Omega(M) = {\rm ker}\, p $.  As the number of rhombi at a given vertex depends on the type of vertex,  
we must deal with the cases individually.   
 
\bigskip

\begin{lemma} {\rm (3-vertex projectives)}   
Suppose $z$ is a 3-vertex. Then we have   
$$0\to b_1A\to e_zA \to \bar{c}_2\bar{b}_3A \to 0$$ 
and  
$$ 0 \to \bar{d}_1\bar{b}_1A  
\to e_zA \to b_3A\to 0.  
$$  
In particular,  $\Omega(\bar{c}_2\bar{b}_3)A \cong b_1A$ and $\Omega(b_3A) \cong  
\bar{d}_1\bar{b}_1A$.  
 \end{lemma}

\bigskip

{\it Proof } Left multiplication by $\bar{c}_2\bar{b}_3$ gives an epimorphism from  
$e_zA$ onto $\bar{c}_2\bar{b}_3A$,   
and since $\bar{c}_2\bar{b}_3b_1=0$ we know that $b_1A$ is contained in the kernel  
of this map.   
The description of the multiplicities of $e_zA$ from the Multiplicity  
Lemma, together with the Arrow Lemma   
show that $\underline{\dim}\, e_zA = \underline{\dim}\, b_1A + \underline{\dim}\,  
\bar{c}_2\bar{b}_3A$, and hence the first sequence is exact.  
Similar arguments yield the second sequence. \hfill $\Box$  
  
\bigskip

\begin{lemma}\label{4-vertex projectives} {\rm (4-vertex projectives)}   
 Suppose $z$ is a 4-vertex. Then there is  
a short exact sequence  
$$0\to b_1A  \to e_zA \to \bar{b}_3A \to 0  
$$  
In particular,  $\Omega(b_1A) \cong  \bar{b}_3A$.  
\end{lemma}

\bigskip

\bigskip

{\it Proof } The arrow $\bar{b}_3$ ends at vertex $z$, and  
we have a surjection $e_zA\to \bar{b}_3A$ given by left multiplication  
with $\bar{b}_3$. Since $\bar{b}_3b_1=0$, we see that $b_1A$ is contained in   
the kernel of this map. Again, the Multiplicity Lemma and   
the Arrow Lemma imply that $b_1A$ must be equal to the kernel. \hfill $\Box$

\bigskip

\begin{lemma}\label{5-vertex projectives} {\rm  (5-vertex projectives) }   
 We have a short exact sequence of $A$-modules  
$$0\to b_1A + b_{2}A \to e_zA \to \bar{b}_4A  \to 0  
$$  
Hence $e_zA$ has a rhombus filtration.  
Moreover, $\Omega(\bar{b}_4A) \cong b_1A+b_2A$.  
\end{lemma}

\bigskip  
  
{\it Proof }   We have a surjective homomorphism  
$e_zA \to \bar{b}_4A$ given by left multiplication. Then $b_1A+b_2A$   
is contained in the kernel. It follows from the Two-arrow Lemma  that   
$b_1A+b_2A$ has a rhombus filtration with three quotients. Furthermore,  
from the Arrow Lemma we know that $\bar{b}_4A$ has a rhombus filtration  
with two quotients. In each case we know the composition factors,  
and it follows that the sequence is exact.  
In particular, $e_zA$ has a rhombus filtration. \hfill $\Box$  
  
\bigskip

\begin{lemma} {\rm (6-vertex projectives) }  
Assume $z$ is  
a 6-vertex. There is a short exact sequence  
$$0\to b_1A+b_2A \to e_zA \to  (\bar{b}_4, \bar{b}_5)A \to 0.  
$$  
The kernel and the cokernel both have a rhombus filtration with  
three quotients.  
\end{lemma}

\bigskip

{\it Proof } There is an obvious   
surjection from $e_zA$ onto $W:= (\bar{b}_4, \bar{b}_5)A$,  
and $b_1A+b_2A$ is contained in the kernel.  
From the Dual two-arrow Lemma it follows that $W$ has a rhombus filtration with three quotients, and  
we also know that $b_1A + b_2A$ has a rhombus filtration with   
three quotients (see the Two-arrow Lemma). The Multiplicity Lemma now shows that  
the sequence is exact. \hfill $\Box$

\bigskip 
 
{\bf Remark } By the same arguments as in the proofs of  
Lemmas~\ref{4-vertex projectives} and Lemmas~\ref{5-vertex projectives}, if $z$ is a 6-vertex, then there also is a short exact sequence  
$$0\to b_1A+b_2A+b_3A \to e_zA \to   \bar{b}_5A \to 0  
$$ 
such that the kernel has a rhombus filtration with four quotients, the cokernel a rhombus  
filtration with one quotient and $\Omega(\bar{b}_5A) \cong 
 b_1A+b_2A+b_3A$.

\section{Truncations}

 We now examine the finite-dimensional Rauzy algebras $A_i=e_iAe_i$ in the  
light of rhombus filtrations.  
Recall that if $z$ is a vertex in $V_i$, then $e_z e_i=e_z$, and furthermore we have that $e_zAe_i$  
is the indecomposable projective (and injective)  
module of  
$A_i$ corresponding to $z$. In particular, it has  simple top and socle.  
  
\medskip  
  
The functor $(-)e_i$ from the module category of right $A$-modules to the module category 
of right  $A_i$-modules is exact. As a result, for  
any vertex $z$ of $Q_i$, the projective  
$e_zA_i$ has a filtration by truncations of rhombus modules, and these are  
easy to write down. Furthermore, if $R$ is a rhombus module  
and $Re_i \neq 0$, then $Re_i$ is indecomposable (this follows from the shape  
of $Q_i$).   
  
\medskip  
  
\subsection{Truncation of arrow modules}  
 
The truncation  of an  arrow module $bA$ when the arrow $b$  
starts or ends at some vertex $z$ in $Q_i$   
has interesting properties.   
If $b$ starts at vertex $z$, then $bAe_i$ is a submodule of $e_zAe_i$, and hence  
$bAe_i$ has a simple socle. Similarly, if $b$ ends at $z$,    
then $bAe_i$ has a simple top. In either case,  
$bAe_i$ is indecomposable.

The module $bA$ has a   
filtration different from the rhombus filtration which gives more precise information about $bAe_i$.

\begin{figure} [ht] 
\centering 
\input{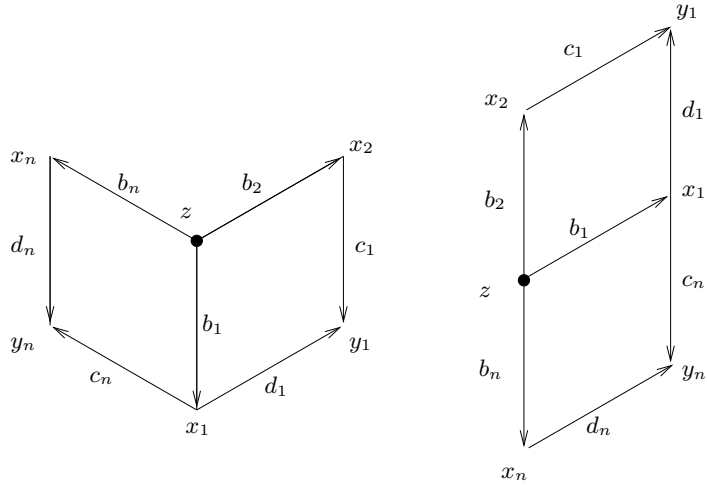} 
\caption{An illustration of adjacent acute and obtuse angles} 
\label{rhomb4a} 
\end{figure}

Let  $b=b_1$ with $b: z \to x_1$. We need to distinguish between two cases determined by the angles at  
$z$ of the two rhombi which have $b$ in common (see figure~\ref{rhomb4a}).  
 
\medskip  
  
(i) either both angles are obtuse, or both angles are acute;  
  
(ii) one of the angles is obtuse, and the other is acute.

\medskip

(i) Consider the first case. Then we have an exact sequence  
$$0 \to \bar{b}bd_1A + \bar{b}bc_nA \to \bar{b}A   
\to  b\bar{b}A \to 0.  \leqno{(1i)}.  
$$  
The kernel is 3-dimensional with a simple socle and top of length two,  
and its composition factors are $S_{y_1},  S_{y_n}$ and $S_{x_1}$.   
There is also an exact sequence   
$$0\to \bar{b}bA \to \bar{b}A \to (\bar{b}_2b\bar{b}, \bar{b}_nb\bar{b})A\to 0. \leqno{(2i)}    
$$  
Here, the cokernel is 3-dimensional  with a simple top and socle  
of length two. Its composition factors are $S_{x_2}, S_{x_n}$ and $S_z$.  
  
\bigskip

Dually, we have exact sequences  
$$0\to b\bar{b}b_nA + b\bar{b}b_2A \to bA \to \bar{b}bA \to 0   \leqno{(3i)}
$$  
with kernel of length three  and composition factors $S_{x_n}, S_{x_2}$ and  
$S_z$; and   
$$0\to b\bar{b}A \to bA \to (\bar{d}_1\bar{b}b, \bar{c}_n\bar{b}b)A\to 0  \leqno{(4i)}  
$$  
with cokernel of length three, and composition factors  
$S_{y_1}, S_{y_n}$ and $S_{x_1}$.

\medskip  
  
(ii) Consider the second case, which can only occur if $z$ is a $5$-vertex or a $4$-vertex.  
 
(a)  Assume $b$ starts at a  
5-vertex $z$, as in the right of figure~\ref{rhomb4a} with $n=5$ and $b=b_1$.  
Then the 
star relation at the vertex $z$ gives  
$$b\bar{b} - b_4\bar{b}_4 = \pm b_3\bar{b}_3.$$ 
Since $\bar{b}_4b=0$ and $\bar{b}_3b=0$,  right multiplication by $b$ shows that $b\bar{b}b=0$. 
 
Similarly we have $\bar{b}b_4=0$ and $\bar{b}b_3=0$, and hence 
$\bar{b}b\bar{b}=0$.  
  
(b)  Assume $b$ starts at a 4-vertex $z$ and a rhombus adjacent to $b$  
has an obtuse angle at $z$. As follows from the general list of possible configurations 
of rhombi as detailed in Table 1 of  section~\ref{horizontal}, the only type of 4-vertex  
which occurs in the Rauzy algebra 
is the completition to a $4$-vertex as in the left of  figure~\ref{rhomb4a}. 
 
Assume $b=b_2$. Then the star relation at $z$ gives  
$$b\bar{b} = \pm b_3\bar{b}_3,$$ 
and  $\bar{b}_3b=0$  implies $b\bar{b}b=0$. 
 
Now we have seen that in case (ii) $b\bar{b}b=0$, and similarly 
$\bar{b}b\bar{b}=0$. Therefore we have the exact sequences  
$$0\to \bar{b}bA \to \bar{b}A \to b\bar{b}A\to 0, \leqno{(1ii)}  
$$  
$$0\to b\bar{b}A \to bA \to \bar{b}bA\to 0.  \leqno{(2ii)} 
$$  
  
\bigskip

\begin{lemma} {\rm (Arrow truncation)}   
Let $z$ be a $4$-vertex.

(a) Suppose $e_{x_1}e_i=0=e_{y_1}e_i= e_{y_4}e_i$ but $e_ze_i\neq 0$.   
Then we have $bAe_i =b \bar{b}Ae_i=b\bar{b}A_i = \bar{b}Ae_i$.   
  
(b) Suppose $e_ze_i=0=e_{x_4}e_i = e_{x_2}e_i$ but $e_{x_1}e_i \neq 0$. Then   
we have $bAe_i = \bar{b}bAe_i =  \bar{b}bA_i=\bar{b}Ae_i$.  
\end{lemma}  
  
Note that $b$ and $\bar{b}$ do not belong to $A_i$, but in part (a) we have $b\bar{b}$ $\in A_i$ and 
in part (b) we have $\bar{b}b$ $\in A_i$.

\bigskip  
  
{\it Proof } (i) Assume first that $z$ is as in case (i). 

For part (a), apply the exact functor  
$(-)e_i$ to the exact sequence (1i). The kernel becomes  
zero by the hypotheses of part (a), and therefore $$\bar{b}Ae_i = b\bar{b}Ae_i = b\bar{b}A_i .$$  
Furthermore, we take the  exact sequence (4i), which has $bA$ as the middle  
term, and apply the functor $(-)e_i$. By the hypotheses of part (a), this time the cokernel becomes
zero. So  we have $bAe_i = b\bar{b}Ae_i =  b\bar{b}A_i$.  
  
Part (b) is proved similarly by using the exact sequences (2i) and (3i).
  
(ii) Now assume $z$ is as in case (ii). 

We apply the exact functor  $(-)e_i$ to the exact sequences (1ii) and (2ii). With the hypotheses of
part (a), we have $\bar{b}bAe_i=0$, and the claim follows. To prove (b), we note that
$b\bar{b}Ae_i=0$ in this case.  \hfill $\Box$

\bigskip  
    
Suppose $0 \neq \beta =  \bar{b}b \in A_i$, with $b$ and $\bar{b}$ arrows in $A$ but not in $A_i$,
then $\beta$ may or may not belong to $ {\rm rad\,}^2 A_i$. This leads to the following definition.
  
\begin{definition}{\rm We call $\beta \in A_i$  a {\it loop} whenever $\beta =\bar{b}b$,   
$\beta \notin {\rm rad\,}^2 A_i $, and $b$ and  $\bar{b}$ are arrows in $A$ but  not  in $A_i$.}   
 \end{definition} 
  
\subsection{Identifying loops in $A_i$}  
  
\bigskip  
To properly understand the structure of $A_i$, we must determine which  
vertices of $Q_i$ admit loops.

\begin{definition} {\rm We call a $(k,n)$-vertex $z$ {\it isolated} if $k = 2$.  
We call an isolated vertex  {\it acute} (resp. {\it obtuse}) if the corresponding angle at $z$ is acute 
(resp. obtuse).}  
\end{definition} 
 
Note that isolated vertices always lie on the boundary 
of the corresponding $\p_i$. 
We use the notation of  section \ref{tiling}, and rhombi are considered to be of type  
$1,2$ or $3$ as indicated in figure~\ref{step0}.

A $\C$-patch denotes a translation of $\p_0$ occuring as a patch in some $\p_i$.  
And a $\C'$-patch denotes a translation of the patch corresponding to $\mathcal{ R}\left(R_1 \right).$  
 
\bigskip 
 
\begin{lemma}\label{Lem}  
If a rhombus $R$ in $\p_i$ occurs as the image under either one or two iterations of  
$\mathcal{ R}$ of a rhombus in some  $\C$-patch or $\C'$-patch, then there is no acute isolated vertex in $R$.   
\end{lemma}  
 
\bigskip 
 
{\it Proof} The function $\mathcal{ R}$ induces a function on patches as well as on individual  
tiles. Since the image under either one or two iterations of  $\mathcal{ R}$ on either a $\C$-patch or a $\C'$-patch has no isolated vertex and since the occurence of a  
patch in $\p_i$ can only transform isolated vertices into unisolated vertices, one obtains the lemma. \hfill $\square$

\begin{lemma}\label{acute}  
For any $i \geq 0$, no acute isolated vertex occurs in $\p_i.$  
\end{lemma}  
  
{\it Proof }        The proposition is clear for $i=0.$  
  
        For $i\geq 1$, a rhombus in $\p_i$  of type 3 can only occur as the image of a rhombus of type  
        1. Hence, any type 3 rhombus  
in $\p_i$ must occur as part of a $\C'$-patch, which can have no acute isolated vertex.

        For $i\geq 1$, a rhombus $R$ of type 2 can only occur as an image of a rhombus of type 1 or 3. If $R$ is the image of  type 1 rhombus, $R$ must occur as part of a $\C'$-patch. If $R$  
occurs as the image of a rhombus $R'$ of type 3 in $\p_{i-1},$ then $R'$  must be part of a $\C$-patch  
or a $\C'$-patch, and so by  lemma~\ref{Lem} $R$ can have no isolated acute vertex.

        For $i\geq 1$, a rhombus $R$ of type 1 can only occur as an image of a rhombus of type 1  
or 2.  If $R$ is the image of a type 1 rhombus, $R$ must occur as part of a $\C'$-patch. If $R$  
occurs as the image of  $R'$ of type 2, by considering the previous case, one sees that  
lemma~\ref{Lem} applies to $R$ as well. \hfill $\square$  
  
\bigskip  
  
\begin{lemma}\label{nkvertice} For any  $(k,n)$-vertex $z$,  
$n-k \leq 3$.  
\end{lemma}  
{\it Proof}  As every vertex is the vertex in a rhombus, there are always at least two edges and at  
most six edges at $z$;  thus, $2 \leq n,k \leq 6$. Therefore we always have $n-k \leq 4$, and we only  
have  
to check that $n-k =4$ does not occur. As $n \geq k$, the only possibility is a $(2,6)$-vertex.  
But a $(2,6)$-vertex  corresponds to an acute isolated vertex in $\p_i$, which cannot occur by  
lemma \ref{acute}. \hfill $\square$

\begin{proposition}  
The  $(k,n)$-vertex $z$ yields a loop in  $A_i$ if and only if $k < n-2$.  
\end{proposition}

{\it Proof } As there are $n$ different  
paths of length two from and to  $z$ in $Q$ and there are two independent  
relations, the dimension of the space $X_z$ is $n-2$.  
  
Then ${\rm rad \,}^2(A_i)\cap X_z$ is spanned by  
$k$ elements and hence has dimension $\leq k$. If $k<n-2$, then this space is  
strictly contained in  $X_z$  and there are loops. If $n-2-k=1$ then there is just one. 

Suppose $k\geq n-2$. One checks easily, case by case, that always ${\rm rad \,}^2(A_i) \cap X_z =
X_z$, and hence there is no loop in this case. \hfill $\square$   
\medskip

By lemma \ref{nkvertice},  $(3,6)$-vertices and   
$(2, 5)$-vertices are the only such vertices. In these cases there is just one loop at such a vertex.   
  

{\rm If $P$ is indecomposable projective and injective, the subquotient 
${\rm rad} P/{\rm soc } P$ is often called the `heart' of $P$. }

\medskip  
  
\begin{proposition}\label{loop}If there is a loop  
at  the vertex $z$, then the projective $P_z$ has decomposable heart, and one of the summands  
is simple and isomorphic to $ S_z$. 
\end{proposition}  
  
\bigskip

{\it Proof }Suppose $z$ is a $(3,6)$-vertex.  
Let $x_1, x_2, x_3$ be the neighbouring vertices of $z$ in $Q_i$ and let $y_1, y_2, b_1, \ldots$ be as before  
(see figure~\ref{star6}).  
Assume that $b_1, b_2, b_3$ are neighbouring arrows  
which are nonzero in $A_i$.  
Then $X_z \cap (\sum_{i=1}^3 b_iA)$ is spanned by $b_i\bar{b}_i$ 
for $i=1, 2, 3$, and they are linearly independent.  
But $X_z$ has dimension four, so we need some element in $X_z$ to generate 
the radical of $e_zA_i$.  
 
We can take as the additional generator the path $\gamma = b_5\bar{b}_5$.  
Then $\gamma b_1= \gamma b_2=\gamma b_3=0$.  
 
Let $M:= \sum_{i=1}^3 b_iA$, then ${\rm rad} P_z = M + \gamma A_i$ 
and $M \cap \gamma A_i = \gamma^2A_i = {\rm soc} P_z$.  
That is, ${\rm rad }P_z/{\rm soc} P_z$ is the direct sum of 
a 1-dimensional module with $M/{\rm soc }P_z$, and the 1-dimensional 
summand is isomorphic to $S_z$.

Similarly, if $z$ is a $(2,5)$-vertex, then the loop gives rise to  
a 1-dimensional direct summand of $rad( P_z)/soc(P_z)$. \hfill $\Box$

\bigskip  
  
\section{ The stable Auslander-Reiten quiver of $A_i$}

The  {\it Auslander-Reiten quiver} $\Gamma(\Lambda)$ for a finite-dimensional symmetric $k$-algebra  $\Lambda$  
is the quiver with vertices the isomorphism classes 
of indecomposable $A$-modules and arrows $[M]\to [N]$ corresponding to irreducible 
maps  $f: M\to N$, where $[M]$ is the class of the module $M$. Thus, by providing generators and some relations,   
$\Gamma(\Lambda)$ may be  thought of as yielding part of a presentation of the  module category 
of $\Lambda$. The irreducible maps are occur in the setting of almost split sequences.  
The almost split sequence ending in $N$ 
is a non-split short exact sequence of $A$-modules 
$$0\to \Omega^2(N)\stackrel{g}\to X\stackrel{h}\to N\to 0 
$$ 
such that any map $\Omega^2(N)\to Y$ which is not a split monomorphism 
factors through $g$. 
Equivalently, any map $Y\to N$ which is not a split epimorphism factors 
through $h$. 
 
There is an irreducible map $M\to N$ if and only if 
$M$ is isomorphic to a direct summand of $X$, and if so then the irreducible map occurs 
as a component of $h$. Dually, there is an irreducible map 
$\Omega^2(N)\to M$ if and only if $M$ is a direct summand of $X$, 
and if so then the irreducible map occurs as a component of $g$; 
see, for example, \cite{ARS} or \cite{B}. 
\medskip 
 
Any indecomposable projective module $P$ of $\Lambda$ with  
simple top $S$ has a unique almost split sequence where 
$P$ occurs, namely 
 
$$0 \to \Omega(S) \to P\oplus {\rm rad} P/{\rm soc} P \to \Omega^{-1}(S)\to 0 
\leqno{(*)} $$ 
 
The  {\it stable Auslander-Reiten quiver} $\Gamma_s(\Lambda)$ of $\Lambda$ is obtained from 
$\Gamma(\Lambda)$ by removing the projective modules and  adjacent arrows. 
Any connected component of $\Gamma_s(\Lambda)$ is, as a graph, isomorphic 
to $\mathbb{Z} T/ \Pi$, where $\Pi$ is some admissible group of automorphisms. The graph 
$\bar{T}$ associated to $T$ is uniquely determined by the component, and 
is called its {\it tree class}, hence  
one wants to identify $\bar{T}$.  
The main method to do this, for symmetric (or self-injective) algebras, uses 
subadditive functions, and for these one needs 
to establish the  existence of suitable $\Omega$-periodic modules, which we herafter 
refer to as {\it periodic modules}. For more detail, see~\cite{B} and~\cite{EHSST}.

\begin{lemma} \label{types} {\rm(\cite{EHSST})} Suppose $\Theta$ is a component  
of the stable Auslander-Reiten quiver of a symmetric algebra. 
If there is some periodic module $W$  such that  
$\underline{\rm Hom}_{\Lambda}(W, X)\neq 0$ for some module $X$ in $\Theta$, 
then the tree class of $\Theta$ is Dynkin, Euclidean, or  
one of the infinite trees $A_{\infty}, A_{\infty}^{\infty}$,  $D_{\infty}$, 
or possibly $B_{\infty}$ or  $C_{\infty}$. 
\end{lemma}

\medskip 
 
With this hypothesis, one takes as (sub)additive function 
the map $$d_W(-):= \dim(\underline{\rm Hom}(\hat{W}, -),$$ where 
$\hat{W}$ is the direct sum of all distinct $\Omega$-translates of $W$.  
 
\medskip

Types $B_{\infty}$ or $C_{\infty}$ cannot occur for algebraically closed fields. 
For group algebras or other symmetric algebras arising naturally where tree classes are known,  
one usually has the tree class of $A_{\infty}$, and one might expect 
this to be the case more generally. Now we focus our attention again on 
the algebras $A_i$ to determine whether suitable (sub)additive functions exist 
and to see if the tree classes of simple modules can be identified.  
 
\medskip

\begin{proposition} Suppose $x$ is a vertex of the quiver of $A_i$, and 
assume that there is a periodic module $W$ such that 
$\underline{\rm Hom}(W, S_x)\neq 0$.  
Let $H_x := {\rm rad}(P_x)/{\rm soc}(P_x)$, and let $\bar{T}$ be the 
tree class of the stable component containing $S_x$.  
If $\bar{T}$ is an infinite tree, we have the following dichotomy. 
 
(a) If $H_x$ is indecomposable, then $\bar{T}\cong A_{\infty}$.  
 
(b) If $H_x$ is decomposable, then the tree class is not $A_{\infty}$.  
\end{proposition}

\bigskip 
 
{\it Proof } 
 (a) We work over $\bar{k}$, the algebraic closure of  $k$. In this case,   
$B_{\infty}$ and $C_{\infty}$ are excluded.  
If $\bar{T}=A_{\infty}$ over  $\bar{k}$, then it follows from  
`Galois descent' (see~\cite{B}) that $\bar{T} = A_{\infty}$ over $k$ as well.

\medskip 
The components of $S_x$ and of $H_x$ have the same tree class. 
This holds since $H_x$ is in the component of $\Omega(S_x)$, as the standard 
sequence shows, and $\Omega$ induces a stable self-equivalence 
of $A_i$, hence it induces a graph isomorphism of the 
stable Auslander-Reiten quiver. 
 
So we consider the tree class of the component of $H_x$ and therefore of $\Omega(S_x)$.  
Recall that the vertices of $\bar{T}$ correspond to the $\Omega^2$-orbits 
on the component. 
Since $H_x$ is indecomposable, the middle term of 
the almost split sequence starting in $\Omega(S_x)$ has indecomposable 
non-projective part. This means that the corresponding vertex of $\bar{T}$ 
must be an end vertex.  
 
Suppose  $\bar{T}\neq A_{\infty}$, then the only other possibility is 
 $D_{\infty}$. Since the vertex of $\bar{T}$ corresponding to $\Omega(S_x)$ 
 is an end vertex, it follows that the vertex of $\bar{T}$ corresponding 
 to $H_x$ is the branch vertex. Therefore, by the definition 
 of $\mathbb{Z}D_{\infty}$, there is an almost split sequence of the 
 form  
$$0 \to \Omega^2(U)  \to H_x \stackrel{g}\to U \to 0 
\leqno{(**)} 
$$ 
where $U\not\cong \Omega^{-1}(S_x)$.

We have ${\rm soc}( H_x )\cong  \oplus S_{x_i}$, where 
$x_i$ ranges over all the neighbours of $x$ in $Q_i$.  
Now, ${\rm soc}( \Omega^2(U))$ is contained in ${\rm soc}( H_x )$. On the 
other hand, we have ${\rm Ext}^1(U, S_x) \neq 0$, and a non-zero element in 
this space exists, namely the 
inverse image of $U$ in $P_x$. But  
$${\rm Ext}^1(U, S_x) \cong \underline{\rm Hom}(\Omega(U), S_x) 
\cong {\rm Hom}(\Omega(U), S_x) 
$$ 
and hence $S_x$ occurs in the top of $\Omega(U)$, which is isomorphic to 
the socle of $\Omega^2(U)$, a contradiction. 
 
\bigskip

(b) Suppose now that $H_x$ is a direct sum. This happens only when 
the quiver of $A_i$ has a loop at vertex $x$.   
In this case we know that $H_x$ has a direct summand 
isomorphic to $S_x$. Then by considering the standard sequence (*) 
we see that there is an irreducible map $\Omega(S_x) \to S_x$.  
This shows that the graph automorphism of the stable Auslander-Reiten-quiver  
fixes the component of $S_x$, but it induces a non-trivial automorphism 
of the tree $\bar{T}$ as the $\Omega^2$-orbits of $S_x$ and of $\Omega(S_x)$ are 
distinct.  
But $A_{\infty}$ does not have a non-trivial automorphism. \hfill $\Box$

\medskip

By Auslanders's theorem, an indecomposable agebra of infinite type does not have a component where the 
tree class is Dynkin.  
 
\medskip 
 
{\bf Conjecture } {\it (a) For any simple module $S_x$ there is a periodic module $W$ with 
$\underline{\rm Hom}(W, S_x) \neq 0$ . 
 
(c) Euclidean components do not occur in the stable Auslander-Reiten-quiver of $A_i$.  
 
(c) If $A_i$ has a loop at the vertex $x$,  
then the tree class of the component of $S_x$ is 
$A_{\infty}^{\infty}$. }

\medskip 
 
To prove this, we would need suitable periodic modules. 
We will next show that we have periodic modules, and this will 
allow us to identify tree classes of some components of simple modules. 
  
\bigskip

\section{Search for periodic modules}

In this section we examine when arrow modules are periodic.  
 
\medskip  
  
\subsection{$\Omega$-Translates of arrow modules} \label{omegatranslates} 
  
For an arrow $b:v \to z$, the arrow module depends on $z$. Below is a list of the different possibilities. 
 
\medskip 
 
{\bf Arrow translates:}\label{arrowtrans} 
\begin{enumerate}   
\item If $b$ is an arrow which ends at a   
5-vertex $z$ (see figure~\ref{star} with $n=5$ and $b=\bar{b}_1$), then we have seen that   
$\Omega(bA) = b_3A + b_4A$.  
 
\item If $b$ ends at a 4-vertex (see figure~\ref{rhomb4a} with $b=\bar{b_1}$), then $\Omega(bA)=b_4A$. (This then determines the  $\Omega$-translates for arrow modules along any path in the quiver  
that goes only between 4-vertices).  
 
\item If the arrows $b_3$ and $b_4$ in (1) bound the same rhombus and both  
end at a $3$-vertex, then   
$\Omega(b_3A + b_4A) = (c_3, -d_3)A$.  

\item If two adjacent arrows $\bar{b}_1$ and $\bar{b}_2$ end at a 5-vertex (see figure~\ref{star} with $n=5$
  with an obtuse angle between $b_4$ and $b_5$ or $b_3$ and $b_4$), then 
$\Omega\left((\bar{b}_1,\bar{b}_2)A\right) = b_4A$.

\item If two adjacent arrows  $\bar{b}_1$ and $\bar{b}_2$  end at a 6-vertex (see figure~\ref{star6}), then  
$\Omega\left((\bar{b}_1,\bar{b}_2)A\right) = b_4A + b_5A$.

\end{enumerate}  
  
\subsection{Horizontal lines}\label{horizontal}  
 
In order to identify periodic modules, we need to introduce some terminology that  
describes the geometry of the $\Omega$-translates of arrow modules. 
 
\medskip 
 
\begin{figure}[ht]  
\centering  
\epsfig{file=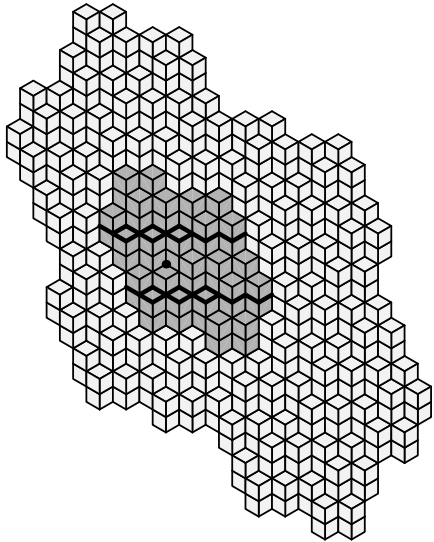, width=6cm}  
\caption{$\p_6$ within $\p_9$ and two horizontal lines corresponding to periodic modules. }  
\label{6in9}  
\end{figure}

\begin{figure}[ht]  
\centering  
\epsfig{file=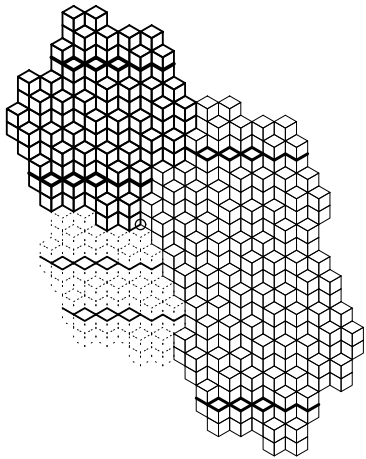, width=5cm}  
\caption{$\p_9$ with relevant  horizontal lines indicated in each $\mathcal{R}^9\left(R_i\right)$ }  
\label{R9}  
\end{figure}

We call an edge in $\p_i$ {\it horizontal} if it is not vertical; that is, if it is parallel to  
$\mathrm{p} \left(\mathbf{e}_1\right)$ or $\mathrm{p} \left(\mathbf{e}_2\right)$. We say that two vertices  
are {\it on the same level} if there is a sequence of horizontal edges joining them. Then we 
refer to a collection of edges between two vertices $v$ and $z$  on the same level as a {\it horizontal line}  
if it contains all horizontal edges joining vertices that are on the same level as and  between $v$ and $z$.  
Figures~\ref{6in9} and~\ref{R9} provide some  
examples of horizontal lines in bold.  Finally, the {\it width} of $Q_i$ at the level of a vertex $z$ is 
the minimum number of horizontal edges joining boundary points on the same level as $z$. 
 
\medskip 
 
This new terminology allows us to describe the results of the previous section in greater detail as follows. 
 
\medskip 
 
\begin{proposition} If $b$ is a horizontal arrow in $Q_i$ which starts at the left or right boundary of $Q_i$  
and  if $k$ is the width of $Q_i$ at the level of $b$, 
then the first $k$  $\Omega$-translates of $bA$ are one of the arrow translates as described in  
section~\ref{arrowtrans}  and lie along a horizontal line. 
  
\end{proposition}  
  
\bigskip  
  
We now study the $\Omega$-translate of a module  
as the horizontal line reaches the boundary of $Q_i$.  
The important question is whether the $\Omega$-translates  
after step $k$ are still small.  
  
\bigskip  
 
\begin{lemma}\label{omegalem}  
Let  $b$ be an arrow in $Q_i$ starting at a $4$-vertex and ending at the $4$-vertex $z$. Suppose that in $A$  
we have $\Omega(bA)= hA$ for some path $h$ of $A$ such that $he_i=0$. 
Then $\Omega(bA_i) \cong  h\bar{h}A_i$ and $\Omega^2(bA_i) \cong  
\bar{b}A_i$.  
\end{lemma}

{\it Proof } 
We then have an exact sequence of $A$-modules 
$$0\to hA\to e_zA \to bA\to 0.  
$$  
The module $bAe_i= bA_i$ is  still an arrow module in the algebra $A_i$.   
The functor $(-)e_i$ is exact, and therefore $e_zAe_i$, the projective indecomposable attached to $z$ in  
$A_i$, is a projective cover of $bA_i$.    
Hence in $A_i$ the $\Omega$-translate of $bA_i$ is equal  
to $hAe_i$. It follows from the Arrow truncation Lemma  that   
$$hAe_i = h\bar{h}Ae_i = \bar{h}Ae_i.  
$$  
Furthermore, we have an exact sequence in $A$  
$$0\to \bar{b}A \to e_zA \to \bar{h}A\to 0.  
$$  
Applying the exact functor $(-)e_i$ gives the projective cover for $\bar{h}Ae_i$, demonstrating 
that $\Omega(\bar{h}Ae_i) = \bar{b}A_i$. \hfill $\Box$  
 
\medskip  
 
Hence, we have the following.

\begin{theorem}\label{perpath}  Suppose the horizontal line $L$ in $Q_i$ 
joins two $4$-vertices on the boundary. Let $b$ be an arrow on $L$ 
starting or ending at a $4$-vertex. Then $bA_i$ is periodic of period   
$2k+2$, where $k$ is the width of the quiver at this level. Furthermore,   
the $\Omega$-translates of $bA_i$  lie along $L$: two modules are as in lemma~\ref{omegalem} 
(corresponding to the ends of $L$), while the other translates are as  
described in section~\ref{arrowtrans}. 
\end{theorem} 
  
\medskip  
  
Thus, we have an entirely combinatorial criterion for periodicity that allows us to identify  
periodic modules in each $A_i$ for sufficiently large $i$.  
  
\medskip  
  
\begin{theorem} For each $i \geq 6,$ $A_i$ contains a periodic module. 
\end{theorem} 
{\it Proof }  
We demonstrate the existence of periodic modules in each algebra $A_i$ for $i \geq 6$ 
directly.  First, we show how to find horizontal lines corresponding to two such modules in  
$A_i$ $6 \leq i \leq 9$, then we  
describe how the construction of $A_{i+1}$ from $A_i$ guarantees the existence of two similar  
horizontal lines for all $i \geq 9$. Finally, by examining the limited number of configurations  
of tiles in the  
Rauzy tiling, we show that these horizontal lines will always yield periodic modules.

Consider the two indicated horizontal lines  in $\p_6$ that end in what can be seen to be 4-vertices by  
examining the  
surrounding tiles in $\p_9$ as illustrated in figure~\ref{6in9}. These horizontal lines correspond to periodic modules. Translations of these horizontal lines can be  
found in $\p_7$ and $\p_8$ joining 4-vertices on the boundary. But in $\p_9$, there are such  
horizontal lines  in \emph{each}  $\mathcal R^9 \left( R_k  
\right)$ joining 4-vertices of the boundary of the corresponding patch  $\mathcal R^9 \left( R_k  
\right)$, as indicated in figure~\ref{R9}.

Since  
$\mathcal R^{10} \left( R_k \right)$ for $k=2,3$ are simply  $\mathcal R^9 \left( R_k \right)$ for  
$k=1,2$, these patches will also contain such horizontal lines joining 4-vertices on the boundary. As $\mathcal R^{10} \left( R_1 \right)$ is  
formed as a combination of translations of the  $\mathcal R^9 \left( R_k  
\right)$, it too will have such horizontal lines. By induction, this will be the case for each subsequent  
$\p_i$. Then for each such $i$, $\p_i$ will have at least two horizontal lines going from one boundary  
vertex to another. (Precisely which patches $\mathcal R^i \left( R_k \right)$ contribute these horizontal lines  
depends on the class of $i \mod 3$).  
  
However, in order to conclude that these horizontal lines correspond to periodic modules, we must know that  
each of these horizontal lines begin and end at a 4-vertex. For this we will need the results of \cite{BV},  
which tells us how many configurations of rhombi of a certain size  the Rauzy tiling contains.  
 Recall the  
lattice $\mathcal{L}$, on which the vertices of the rhombi lie, has basis $\left\{\mathbf{b}_1 = \mathrm{p} \left(\mathbf{e}_2\right)\,,  
  \mathbf{b}_2 =  \mathrm{p} \left(\mathbf{e}_3\right) \right\}.$ When we consider the rhombi of the three types  to have bases as indicated in figure~\ref{bases},  
each vertex of  $\mathcal{L}$ is the base for exactly one rhombus~\cite{BV}.

\begin{figure}[h!]  
\centering  
\epsfig{file=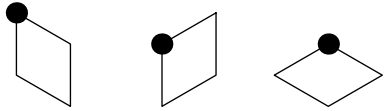, width=6cm}  
\caption{Each $R_k$ with its base indicated}  
\label{bases}  
\end{figure}

Then there are $mn+m+n$ different configurations of rhombi based in $\left(m,n\right)$-sections of   
$\mathcal{L}$ composed of vertices of the form  
\[  
\left\{(k+i)\,\mathbf{b}_1 + (\ell +j)\, \mathbf{b}_2 : i=1,\dots,m \, ; j=1, \dots,n \right\}  
\]  
for some $k$ and $\ell$~\cite{BV}.

\begin{figure}[h!]  
\centering  
\epsfig{file=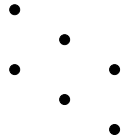, width=1.5cm}  
\caption{Example of a $(3,2)$-section of $\mathcal{L}$}  
\label{bases}  
\end{figure}

Below we include a table of these configurations for select $\left(m,n\right)$, where the entry  
$a_{ij}$  
of the matrix representing a configuration corresponds to the type of rhombus (1, 2 or 3)  
based at $(k+i)\,\mathbf{b}_1 + (\ell +(n-j))\, \mathbf{b}_2$.  
 
\vspace{1cm}

\begin{tabularx}{\linewidth}%
     {|>{\setlength{\hsize}{.3\hsize}}X|%
       >{\setlength{\hsize}{.3\hsize}}X|%
       >{\setlength{\hsize}{2.4\hsize}}X|} 
\hline 
$(m,n)$ & $P(m,n)$ &  List of possible configurations of rhombi \\  \hline\hline 
$(1,1)$ & $3$        &   $(1)$, $(2)$, $(3)$  \\ \hline 
$(1,2)$ & $5$ &  $\left(\begin{array}{c} 1 \\ 1 \end{array}\right)$,  
$\left(\begin{array}{c} 2 \\ 2 \end{array}\right)$,  
$\left(\begin{array}{c} 1 \\ 3 \end{array}\right)$,  
$\left(\begin{array}{c} 2 \\ 1 \end{array}\right)$,  
$\left(\begin{array}{c} 3 \\ 2 \end{array}\right)$ \\ \hline  
$(2,1)$ & $5$ & $(1,1)$, $(1,2)$, $(3,1)$, $(2,1)$, $(2,3)$   \\ \hline 
$(2,2)$ & $8$ &  $\left(\begin{array}{cc} 1& 2 \\ 1 & 2 \end{array}\right)$, 
$\left(\begin{array}{cc} 2& 1 \\ 2 & 3 \end{array}\right)$, 
$\left(\begin{array}{cc} 1& 2 \\ 1 & 1 \end{array}\right)$, 
$\left(\begin{array}{cc} 1& 2 \\ 3 & 1 \end{array}\right)$, 
$\left(\begin{array}{cc} 1& 1 \\ 3 & 1 \end{array}\right)$, 
$\left(\begin{array}{cc} 2& 3 \\ 1 & 2 \end{array}\right)$, 
$\left(\begin{array}{cc} 3& 1 \\ 2 & 3 \end{array}\right)$, 
$\left(\begin{array}{cc} 3& 1 \\ 2 & 1 \end{array}\right)$ 
\\ \hline 
$(3,2)$ & $11$ &  $\left(\begin{array}{ccc} 1 & 2 & 1 \\ 1 & 2 & 3 \end{array}\right)$, 
$\left(\begin{array}{ccc} 1 & 1  & 2 \\  3 & 1 & 2 \end{array}\right)$, 
$\left(\begin{array}{ccc} 1 & 2 & 3 \\ 1 & 1 & 2 \end{array}\right)$, 
$\left(\begin{array}{ccc} 1 & 2 & 3 \\ 3 & 1 & 2 \end{array}\right)$, 
$\left(\begin{array}{ccc} 2 & 1 & 1 \\ 2 & 3 & 1 \end{array}\right)$, 
$\left(\begin{array}{ccc} 2 & 1 & 2  \\ 2 & 3 & 1 \end{array}\right)$, 
$\left(\begin{array}{ccc} 2 & 3 & 1 \\ 1 & 2 & 1 \end{array}\right)$, 
$\left(\begin{array}{ccc} 2 & 3 & 1 \\ 1 & 2 & 3 \end{array}\right)$, 
$\left(\begin{array}{ccc} 3 & 1 & 2 \\ 2 & 1 & 1 \end{array}\right)$, 
$\left(\begin{array}{ccc} 3 & 1 & 2 \\ 2 & 1 & 2 \end{array}\right)$, 
$\left(\begin{array}{ccc} 3 & 1 & 2 \\ 2 & 3 & 1\end{array}\right)$ 
 \\ \hline 
\end{tabularx} 
 
\begin{center} Table 1 
\end{center}

\bigskip

Each of the horizontal lines under consideration ends on the right with a   
$\left(%
\begin{array}{rrr}  
   
  3  \\  
  2 \\  
\end{array}%
\right)$  
configuration above, and so consideration of the possible $\left(2,2\right)$-configurations and  
geometry forces this to be completed as  
$\left(%
\begin{array}{rrr}  
   
  3 & 1 \\  
  2 & 1 \\  
\end{array}%
\right),$  
corresponding to the desired 4-vertex. On the left, the horizontal lines end cutting through a  
$ \left(%
\begin{array}{rrr}  
   
  1 & 2  \\  
  1 & 1 \\  
\end{array}%
\right)$  
configuration, and so consideration of the possible $\left(3,2\right)$-configurations shows  
that it must be completed on the left as  
$ \left(%
\begin{array}{rrr}  
   
 3 &  1 & 2  \\  
 2 &  1 & 1 \\  
\end{array}%
\right).$  
This means that the configuration on the left end of the horizontal line containing  
the type 1 rhombus based  
at the left endpoint has the form  
$ \left(%
\begin{array}{rrr}  
   
 2 &  1 \\  
 x &  y \\  
\end{array}%
\right),$  
which can only be completed as  
$ \left(%
\begin{array}{rrr}  
   
 2 &  1 \\  
 2 &  3 \\  
\end{array}%
\right)$, yielding the desired 4-vertex. \hfill $\square$  
 
\bigskip 
 
\bigskip

We call a horizontal line {\it periodic} if it corresponds to a periodic orbit. In what follows, $L$ denotes 
a periodic line in $Q_i$. 
\begin{definition}{\rm  
We say that an $A_i$-module $W$ is {\it along $L$} if it has  
arrow translates as described in~\ref{arrowtrans} or if it is of the form 
$W= (b\bar{b})A_i$, where $b$ starts at $L$ but does not 
belong to $A_i$ and $b\bar{b}$ is  in $A_i$. } 
\end{definition} 
 
By theorem~\ref{perpath}, the modules along $L$ form one $\Omega$-orbit. 
 
\medskip

\begin{corollary}  Suppose $M$ is an indecomposable non-projective 
$A_i$-module  and suppose for some $W$ along  
the line $L$ we have 
$\underline{\rm Hom}(W, \Omega^t(M)) \neq 0$ for some  
$t\geq 0$. Then the tree class of the component  
of $M$ belongs to the list as in lemma~\ref{types}.  
\end{corollary} 
\medskip 
 
We can completely answer for which simple modules this applies. 
 
\medskip 
 
\begin{theorem} Suppose $S$ is a simple module. 
Then there is a periodic $W$ along a periodic line $L$ 
such that $\underline{\rm Hom}(W, \Omega^t(S))\neq 0$ for 
some $t$ if and only if $S=S_x$ and $x$ is a vertex on $L$. 
\end{theorem}

{\it Proof }  
(1) Suppose $S=S_x$, where $x$ is a 4-vertex or 5-vertex on  
a periodic line $L$. Then $\Omega(S)$ contains $W=bA_i$ for $b$ an arrow starting at 
$x$, hence $W$ is along $L$. An inclusion map does not factor through a projective module, and  
hence  
$$0\neq \underline{\rm Hom}(W, \Omega(S)). 
$$ 
 
(2) Suppose  $S=S_x$ where $x$ is a 6-vertex on a periodic line $L$. Then 
$\Omega(S)$ contains $W= hA_i+kA_i$, where $h, k$ are arrows starting at $x$ and  
$W$ is along  $L$. 
Then similarly $\underline{\rm Hom}(W, \Omega(S))\neq 0$.  
 
(3) Suppose $S=S_x$ where $x$ is a  
3-vertex which occurs along a periodic line $L$. 
Then $S$ is a top composition factor of $W= hA_i+kA_i$, 
where $W$ is along $L$ and where either $h$ or $k$ ends at $x$. 
The epimorphism from $W$  onto $S$ 
does not factor through a projective module. Hence, 
as before, $\underline{\rm Hom}(W, S)\neq 0$.

\medskip 
 
For the converse, suppose that $S$ is a simple module and that there is some $W$ along 
$L$ such that $\underline{\rm Hom}(W, \Omega^t(S))\neq 0$ 
for some $t$. This space is isomorphic to  
$$\underline{\rm Hom}(\Omega^{-t}(W), S) 
\cong {\rm Hom}(\Omega^{-t}(W), S) 
$$ 
where the first isomorphism is by dimension shift, and the second holds 
because $S$ is a simple module of a symmetric algebra and 
$\Omega^{-t}(W)$ does not have  
a projective summand.  
Now, $\Omega^{-t}(W)$ also is a module along $L$. Since there 
is a non-zero homomorphism from this module onto $S$ it follows that 
$S$ is a top composition factor. But all modules along $L$ 
have top composition factors $S_x$ with $x$ on $L$. Hence $S=S_x$  
as stated. \hfill $\Box$

\bibliographystyle{plain}

\begin{thebibliography}{}    
\bibitem{ABS}  P. Arnoux, V. Berth\'{e}, A. Siegel, Two-dimensional iterated morphisms and discrete planes.  
Theoret. Comput. Sci. \textbf{ 319}  (2004),  no. 1-3, 145--176.  
 Bull. Belg. Math. Soc. Simon Stevin \textbf{ 8}  (2001),  no.~2, 181--207.  
\bibitem{ARS} M. Auslander, I. Reiten, S. Smal\o, Representation theory of Artin algebras.  
Cambridge Studies in Advanced Mathematics, \textbf{36}. Cambridge University Press, Cambridge, 1997 
\bibitem{B}  D. J. Benson, Representations and cohomology. I.  
Cambridge Studies in Advanced Mathematics, \textbf{30}. Cambridge University Press, Cambridge, 1998. 
\bibitem{BV} V. Berthé, L. Vuillon, Tilings and rotations on the torus: a two-dimensional generalization of  
Sturmian sequences.  Discrete Math. \textbf{223}  (2000),  no. 1-3, 27--53.  
\bibitem{CT} J. Chuang, W. Turner, Cubist algebras, preprint.  
\bibitem{F}  N. P. Fogg,  Substitutions in dynamics, arithmetics and combinatorics. Edited by V. Berthé,   
S. Ferenczi, C. Mauduit and A. Siegel.  
Lecture Notes in Mathematics, \textbf{1794}. Springer-Verlag, Berlin, 2002.  
\bibitem{EHSST} K. Erdmann, M. Holloway,  N. Snashall,  \O. Solberg,  R. Taillefer, 
Support varieties for selfinjective algebras.  $K$-Theory  \textbf{33}  (2004),  no. 1, 67--87. 
\bibitem{EM} K. Erdmann, S. Martin, Quiver and relations for the principal $p$-block of $\Sigma\sb {2p}$.   
J. London Math. Soc. (2)  \textbf{49}  (1994),  no. 3, 442--462. 
\bibitem{G} P. Gabriel, Auslander-Reiten sequences and representation-finite algebras.  
Representation theory I, Ottawa 1979. Lecture Notes in Mathematics \textbf{831}, Springer-Verlag,  
Berlin/New York 1980. 
\bibitem{P}  M. Peach, Rhombal Algebras, PhD, University of Bristol, August 2004.  
\bibitem{R}  Rauzy, G. Nombres alg\`{e}briques et substitutions.  Bull. Soc. Math. France  \textbf{110}    
(1982), no.~2, 147--178.  
\bibitem{Ri} C. M. Ringel, Rhombal algbras, report, http://www.math.uni-bielefeld.de/~ringel/lectures.html. 
\bibitem{T} W. Turner,  On seven families of algebras, preprint. 
\end{thebibliography}

\end{document}